\documentclass[a4paper,11pt]{article}
\usepackage{amsmath,amssymb, amsthm}

\usepackage{graphicx}
\usepackage{color}
\usepackage{geometry}
\usepackage{psfrag}
\usepackage{epsfig}

\usepackage{a4}
\usepackage{amssymb,latexsym}
\usepackage[mathscr]{eucal}
\usepackage{cite}
\usepackage{url}

\DeclareFontFamily{OT1}{rsfs}{} \DeclareFontShape{OT1}{rsfs}{m}{n}{
<-7> rsfs5 <7-10> rsfs7 <10-> rsfs10}{}
\DeclareMathAlphabet{\mycal}{OT1}{rsfs}{m}{n}

\newcommand{\be}{\begin{equation}}
\newcommand{\ee}{\end{equation}}

\newcommand{\bea}{\begin{eqnarray}}
\newcommand{\beaa}{\begin{eqnarray*}}
\newcommand{\bean}{\begin{eqnarray}\nonumber}

\newcommand{\eea}{\end{eqnarray}}
\newcommand{\eeaa}{\end{eqnarray*}}

\newcommand{\bel}[1]{\begin{equation}\label{#1}}
\newcommand{\beal}[1]{\begin{eqnarray}\label{#1}}
\newcommand{\beadl}[1]{\begin{deqarr}\label{#1}}
\newcommand{\eeadl}[1]{\arrlabel{#1}\end{deqarr}}
\newcommand{\eeal}[1]{\label{#1}\end{eqnarray}}

\theoremstyle{plain}
\newtheorem{theo}{Th�or�me}[section]

\newtheorem{lemm}[theo]{Lemme}
\newtheorem{Lemma}[theo]{Lemma}
\newtheorem{Theorem}[theo]{Theorem}

\newtheorem{conj}[theo]{Conjecture}

\newtheorem{engdef}[theo]{Definition}

\newtheorem{engrk}[theo]{Remark}

\def \R {\mathbb{R}}
\def \Hyp{\mathbb{H}}

\newcommand{\mcE}{{\mathcal E}}

\newcommand{\mcK}{{\mathcal K}}

\newcommand{\mcS}{{\mathcal S}}
\newcommand{\mcF}{{\mathcal F}}

\def \Sphere{\mathbb{S}}

\providecommand{\bysame}{\leavevmode\hbox to3em{\hrulefill}\thinspace}
\providecommand{\MR}{\relax\ifhmode\unskip\space\fi MR }

\providecommand{\href}[2]{#2}

\title{\textbf{\sc{A family of asymptotically hyperbolic manifolds with arbitrary energy-momentum vectors}}}
\author{Julien Cortier\footnote{Max-Planck-Institut f\"ur Gravitationsphysik
(Albert-Einstein-Institut), Am M\"uhlenberg 1, 14476 Golm, Germany, \texttt{jcortier@aei.mpg.de}.}}
\date{}
\begin{document}

\maketitle

\begin{abstract}
A family of non-radial solutions of the Yamabe equation, with reference the hyperbolic space, is constructed using power series.
As a result, we obtain a family of asymptotically hyperbolic metrics, with spherical conformal infinity, 
with scalar curvature greater than $-n(n-1)$, but which are \emph{a priori} not complete. 
Moreover, any vector of $\R^{n+1}$ is performed by an energy-momentun vector of one suitable metric of this family.
They can in particular provide counter-examples to the \emph{positive energy-momentum
theorem} when one removes the completeness assumption.
\end{abstract}

\tableofcontents

\section{Introduction}

For the last decade, the study of asymptotically hyperbolic manifolds has enjoyed a significant effort from the mathematical GR community
as well as from a pure geometric point of view. It is generally motivated by the same questions that arise in the study of asymptotically 
flat manifolds (or asymptotically Minkowskian initial data). Indeed, hyperbolic manifolds, or anti-de Sitter
initial data appear as a ground state in general relativity with negative cosmological constant, and one can define mass-type 
geometric invariants for initial data that approach these models at infinity, and conjecture some positivity results in the spirit of the 
\emph{positive mass theorem}. 

Significant results with the hyperbolic background have been obtained by X. Wang~\cite{Wan01} and by P.T. Chru\'sciel and M. Herzlich~\cite{CH03}.
They both show that a set of global invariants can be defined for asymptotically hyperbolic manifolds, 
forming the \emph{energy-momentum vector}, although this has been done in a more general way in~\cite{CH03}.
They further establish a \emph{positive energy-momentum theorem}, which, like the asymptotically flat positive mass theorem,
ensures that a positivity property of these invariants is satisfied provided the manifold is complete and fulfills the 
\emph{dominant energy condition}. This latter condition translates here as a negative lower bound on the scalar curvature, which is precisely
the value of the scalar curvature for the model considered in each case, $0$ for $\R^n$ and $-n(n-1)$ for $\Hyp^n$.

The aim of the present work is to construct examples of asymptotically hyperbolic metrics which satisfy the above assumptions excepted 
the completeness, and which violate the positivity result of the positive energy-momentum theorem. 

Such results aim at giving a better knowledge of the behaviour of these invariants and the role played by each assumption of the 
positive energy-momentum theorem. Indeed, this theorem has not yet been proven in the most general and satisfying case.
The results obtained by X. Wang and P.T. Chru\'sciel-M. Herzlich indeed require a further assumption of topological nature 
(namely the existence of a spin structure), whereas, in
an attempt to remove this further assumption, L. Andersson, M. Cai and G. Galloway in~\cite{ACG08} needed another
extra condition on the structure at infinity of the manifold.

In the asymptotically flat case, it is straightforward to exhibit (non-complete) Riemannian manifolds of non-negative scalar curvature 
and of negative mass, by considering the exterior of a Schwarzschild manifold of mass $m <0$. In fact, the mass invariant 
of an asymptotically flat manifold (see R. Bartnik,~\cite{Bar86}) is a number, coincides with the parameter $m$ for Schwarzschild manifolds, 
and thus can take any real value. 

However, when the model at infinity is the (complete) hyperbolic space $\mathbb H ^n$, hence with conformal infinity $\Sphere ^{n-1}$,
the set of mass-type invariants is in fact more complicated, since the natural mass-type object to consider 
for an asymptotically hyperbolic manifold is rather the ``energy-momentum'' vector, made of $n+1$ components.

The \emph{positive energy-momentum conjecture} states that this vector has to be timelike and future-directed if the manifold
is moreover complete and has scalar curvature greater than $-n(n-1)$ (which is the scalar curvature of $\mathbb{H}^n$). 

As for the asymptotically flat case, one can exhibit metrics that satisfy these assumptions except the completeness and 
that violate the energy-momentum conjecture. Well-known examples of such metrics are the 
Kottler-Schwarzschild-anti de Sitter metrics with negative mass parameter $m$, given by the expression in local coordinates
$$b_m = \frac{dr^2}{1-\frac{2m}{r^{n-2}} + r^2} + r^2 \sigma_{n-1}\;,
$$
where $\sigma_{n-1}$ is the standard metric on the unit sphere $\Sphere^{n-1}$. Note that we obtain in these coordinates the expression
of the hyperbolic metric for $m=0$.
Their energy-momentum vector is then timelike, past-directed and 
we can moreover obtain any timelike, future or past-directed energy-momentum vector by considering the family of 
\emph{boosted Kottler metrics}, see e.g.~\cite{CD09}. 

But no example of such a metric, which satisfies the assumptions of the positive energy-momentum theorem except the completeness, 
and which has a spacelike or null energy-momentum vector, has been exhibited so far. 
A work~\cite{ST07} from Y. Shi and L-F. Tam shows the existence, in the 3-dimensional case, of A.H metrics of constant scalar curvature $-6$
on the unit ball of $\R^3$,
with energy-momentum vectors whose Minkowskian norm is positive, prescribed up to a small error, 
and with further properties on the existence of horizons.
But this result holds only for timelike, future-directed energy-momentum vectors. 

The aim of the present work is to fill this gap by proving the main result of this paper:
\begin{Theorem}\label{mainthm}
 Let $\mathbf{p}$ be a vector of the Minkowski space $\mathbb R^{1,n}$. Then, there exists a compact set $K \subset \mathbb R^n$ and
a Riemannian, asymptotically hyperbolic metric $g$ defined on $\mathbb R^n\setminus K$, with scalar curvature $R_g \geq -n(n-1)$, such
that the energy-momentum vector of $g$ is well defined and coincides with $\mathbf{p}$.
\end{Theorem}
In fact, to prove this theorem, we construct an explicit family of metrics which have other remarkable properties such as being conformally flat
and having a constant scalar curvature $-n(n-1)$. 
Such a family may then be useful when one wishes to approximate any asymptotically hyperbolic metric of general energy-momentum vector 
by one from this family. This approximation ``near infinity'' is made possible from the gluing method introduced 
by J. Corvino and R.M. Schoen~\cite{Cor00,CS06, CD03}, and it may be used as a first step in a future attempt
to prove the positive energy-momentum theorem, simplifying the structure at infinity. In the same spirit, for the asymptotically
flat case, one can start the Schoen-Yau's proof of the positive mass theorem
by gluing an asymptotically flat manifold of negative mass with a Schwarzschild manifold of (arbitrarily close) negative mass parameter, 
for example using Corvino's gluing result~\cite{Cor00} , or from a direct calculation such as in~\cite{Smontecatini}. 
The result of this is that one can consider only metrics having nicely behaved asymptotics without loss of generality.

\medskip

The present paper is organized as follows:

In section~\ref{sectionAH}, we give the specific material to work with asymptotically hyperbolic manifolds, in particular their various 
definitions, the definition of the energy-momentum vector, and a quick retrospective of the current knowledge 
on the positive energy-momentum theorem.

In section~\ref{sectionthmexistence}, we establish the Theorem~\ref{mainthm}, which is the main of this work. 
We actually look for a family of conformally hyperbolic metrics, 
whose conformal factor satisfies the Yamabe equation such that the scalar curvature is constant and equal to $-n(n-1)$.
We seek an expression for this conformal factor as a series in $1/r$ --where $r$ is the standard radial coordinate-- whose coefficients are 
non-constant functions defined on the $(n-1)$-unit sphere. 

In section~\ref{appli} at last, we derive some applications, one concerning a discussion around a ``positive energy-momentum 
theorem with boundary'' proved by P.T. Chru\'sciel and M. Herzlich in~\cite{CH03}, 
and the other one showing that the family of metrics constructed in the preceeding section
may well be used as a family of models at infinity in the process of gluing, similarly to what is done by P.T. Chru\'sciel and 
E. Delay in~\cite{CD09}.
One can therefore approximate an asymptotically hyperbolic metric $g$ of constant scalar curvature by one with much nicer properties
at infinity, and whose energy-momentum vector is arbitrarily close to the one of $g$.

\paragraph{Acknowledgements}
The author is grateful to the Albert-Einstein-Institut (Max-Planck-Institut f\"ur Gravitationsphysik, Potsdam, Germany) for financial
support and wishes to thank Piotr T. Chru\'sciel, Erwann Delay and Marc Herzlich for 
their guidance during the preparation of this work as well as Romain Gicquaud for useful comments.

\section{Asymptotically hyperbolic manifolds and energy-momentum}
\label{sectionAH}

\subsection{Asymptotically hyperbolic metrics}

For $n \geq 3$, let $(M, g)$ be a smooth $n$-dimensional Riemannian manifold with boundary $\partial M$. 
We will consider such manifolds $M$ which are furthermore non-compact and which have an asymptotic end 
in which $g$ tends to the hyperbolic metric $b$.

These manifolds appear naturally in general relativity as the initial data of asymptotically anti-de Sitter spacetimes.

In the sequel, we give accurate definitions, following the standard notations (see~\cite{theseromain} for example) :

\begin{engdef}
A manifold $(M,g)$ is said to be \emph{conformally compact} if there exists:
 \begin{itemize}
 \item a smooth compact manifold $\bar{M}$, with interior $\mathring{M}$ and with boundary $\partial \bar M$ such that $\partial \bar M$ 
is the union of $\partial M$ and $\partial_{\infty} M$, with $M =  \mathring M \cup \partial M$;

 \item a \emph{defining function}, i.e. a smooth function $\rho : \bar M \rightarrow \mathbb R _+$ such that
  $\rho ^{-1}(0) = \partial_{\infty} M$, the 1-form $d\rho$ does not vanish on $\partial_{\infty} \bar M$, and the metric tensor
   $\bar g := \rho^2 g$ of $M$ extends to a smooth metric on $\bar M$.
  \end{itemize}
\end{engdef}
The boundary component $\partial_{\infty} M$ will often be referred to as the \emph{boundary at infinity}, or as the 
\emph{conformal boundary}, whereas $\partial M$ will be referred to as the \emph{inner boundary}.

\begin{engdef} A conformally compact manifold $(M,g)$ is \emph{asymptotically hyperbolic} (A.H.) if one has the further condition
$|d\rho|_{\bar g} = 1$ on $\partial_{\infty} M$.
\end{engdef}

\begin{engrk}

\begin{itemize}
\item With respect to the above notations, one has in particular $\sec _g \rightarrow -1$ as $\rho \rightarrow 0$: 
the sectional curvature converges to $-1$ near the boundary at infinity. This motivates the ``asymptotically hyperbolic''
terminology.

\item Under these assumptions, there exists a unique defining function $\rho$ and a neighborhood $U$ of $\partial_{\infty} M$ of the form 
$(0,\varepsilon] \times \partial_{\infty} M$ such that, under a diffeomorphism $\Phi$ between $U$ and the complement $M_{ext}$
of a compact region of $M$, the metric takes the form
$$\Phi^* g = \frac{d\rho^2 + h_{\rho}}{(\sinh \rho)^2}
$$
on $U$, where $(h_{\rho})$ is a family of metrics on a $(n-1)$-dimensional compact boundaryless manifold $N$, 
smooth with respect to $\rho$, such that the limit as
$\rho \rightarrow 0$ is a metric $h_0$ on $N$ with constant scalar curvature $R_h = k(n-1)(n-2)$, $k \in \{-1,0,1\}$.

\item It is common (and in fact convenient) to express these definitions using a ``radial'' coordinate $r$ such that
$$\frac1r = \sinh \rho\;.
$$
In these coordinates, the hyperbolic metric $b$ takes the form
\bel{b} b = \frac{dr^2}{1+r^2} + r^2 \sigma_{n-1}\;,
\ee
where $\sigma_{n-1}$ is the standard metric on the unit $(n-1)$-sphere,
whereas the definition for an asymptotically hyperbolic metric would write
\bel{Phig}\Phi^* g = \frac{dr^2}{r^2 + k} + r^2 h_r\;,
\ee
where $\Phi : [R,+\infty) \times N \rightarrow M_{ext}$ is a diffeomorphism and
where $(h_r)_r$ is a family of metrics on $\partial_{\infty} M$, smooth with respect to $r$, such that the limit as
$r \rightarrow +\infty$ is a metric $h_0$ on $\partial_{\infty} M$ with constant scalar curvature $R_h = k(n-1)(n-2)$, $k \in \{-1,0,1\}$.

\end{itemize}
\end{engrk}

There exists a more restrictive notion of asymptotic hyperbolicity, such as presented in \cite{Wan01,ACG08}:
\begin{engdef}\label{sah} An asymptotically hyperbolic manifold $(M,g)$ is \emph{strongly asymptotically hyperbolic} 
(\emph{S.A.H.} for short) if:

\begin{itemize}
\item the conformal infinity $\partial_{\infty} M$ is the $(n-1)$-dimensional unit sphere $\mathbb S^{n-1}$, equipped with its standard metric
 $h_0$;
\item if $\rho$ is a defining function, defined on a neighborhood $U$ of the conformal infinity,
and such that one can write $g = \frac{d\rho^2 + h_{\rho}}{(\sinh \rho)^2}$, 
one then has the asymptotic expansion as $\rho \rightarrow 0$:
    $$h_{\rho} = h_0 + \frac{\rho^n}{n}h + O(\rho^{n+1})\;,
    $$
    where the terms of the expansion can be differentiated twice.
\end{itemize}
\end{engdef}

The tensor $h$ which appears in the above definition is a rank 2 symmetric tensor on $\mathbb S ^{n-1}$,
called the \emph{mass-aspect tensor}.

\subsection{Mass and energy-momentum of asymptotically hyperbolic manifolds}

We recall here how arise the mass and the energy-momentum invariants for an asymptotically hyperbolic (A.H.) manifold.
As for the mass in an asymptotically flat Riemannian manifold, these invariants appear from an Hamiltonian formulation 
of general relativity, as described in~\cite{ADM61,BI04}, in particular they are related to symmetries of the background $(M_0,b)$ 
which is considered (either the Euclidean or the hyperbolic space). 
Indeed, let $\mcK_0$ be the space of the \emph{static Killing initial data} (or \emph{static KIDs}) as defined by P.T. Chru\'sciel and 
R. Beig in~\cite{BC97}, and more precisely for our particular case at the beginning of~\cite{CH03}.
This space of static KIDs here consists in the set of functions $f \in C^{\infty}(M_0)$ such that the spacetime metrics 
$-f^2 dt^2 + b$ are \emph{static} solutions
of the vacuum Einstein equations. The word ``static'' means that the orbits of the Killing vector $\partial_t$ are timelike
and orthogonal to totally geodesic hypersurfaces (here the level sets of the $t$-coordinate function).

Note that in such a spacetime, the vector $\partial_t$ takes the form
$\partial_t = f \nu$, where $\nu$ is the unit (timelike future-directed) vector normal to $M_0$ in the spacetime.  
Then, an important result from V. Moncrief~\cite{Mon76} relates the space of static KIDs to the cokernel of the linearized 
scalar curvature operator evaluated at the background metric $b$ (or for general KIDs to the cokernel of the linearized constraint operator
linearized at the background initial data considered).
To apply this here, we start writing the linearized scalar curvature operator at $b$, tested against a $(0,2)$-symmetric tensor $h$:
$$L_b h := D R(b) h 
$$
and its formal $L^2$-adjoint, against a function $f \in C^{\infty}(M_0)$ :
$$
L_{b}^*f = -(\Delta_b f) b + \mathrm{Hess}_b f - f
\mathrm{Ric}_b \;,
$$
where $\Delta_b$, $\mathrm{Hess}_b$ and $\mathrm{Ric}_b$ are respectively the Laplace-Beltrami operator,
the Hessian and the Ricci curvature tensor of $b$.
Following the notations of~\cite{CH03,thesebenoit}, let $g$ an asymptotically hyperbolic metric in the sense defined above,
and a diffeomorphism at infinity $\Phi$ such that $\Phi^* g$ is asymptote to $b$ (no matter of the topology 
of the conformal boundary). From Moncrief's result, one then has
\bel{mcK0}\mcK_0 = \ker L_b^*\;.
\ee

The mass integrals (see~\cite{CH03,thesebenoit}) then appear for each $V \in \mcK_0$, as limits of flux integrals
\bel{massintegral}H_{\Phi}(V) = \lim_{R \rightarrow +\infty} \int_{r=R} \mathbb U(V,\Phi^* g - b)(\nu_r)dS_r\;, 
\ee
where $\nu_r$ is the unit normal vector to the hypersurface $\{r=R\}$, where $dS_r$ is the induced measure on this hypersurface from $b$ and
where the integrand term reads, for $e = \Phi^* g - b$, as (see~\cite{thesebenoit})
\bel{U(V,e)}\mathbb U(V,e)  = V\left(\mathrm{div}\;e - d(\mathrm{tr}\;e)\right) - \iota_{\nabla V}e + (\mathrm{tr}\;e)dV\;,
\ee
where the divergence and the trace are computed relatively to $b$.

The integrals~(\ref{massintegral}) depend in general on the chosen chart at infinity $\Phi$.
However, in~\cite{CH03,thesebenoit}, it is proven that, for $k=0$ or $-1$, the space $\mcK_0$ is 1-dimensional
and that given an non-trivial element $V_{(0)}$ of $\mcK_0$, the quantity $m = H_{\Phi}(V_{(0)})$
does not depend on $\Phi$   
provided suitable conditions on the decay of $\Phi^* g - b$ at infinity hold (and provided a volume normalization for $k=0$). 
The quantity $m$ is referred to as the \emph{mass} of an A.H. manifold with conformal boundary of negative ($k=-1$) or zero ($k=0$)
scalar curvature. 

The case $k=1$ (conformal boundary of positive Yamabe type) is slightly more complicated. We will only consider the particular case
where $\partial_{\infty}M$ is the $(n-1)$-sphere $\Sphere ^{n-1}$, with $h_0 = \sigma_{n-1}$.
For quotients of $\Sphere ^{n-1}$ by a discrete subgroup $\Gamma$ of its group of isometries, see the discussion in~\cite{CH03}.

In this case of a spherical conformal boundary, $\mcK_0$ is $(n+1)$-dimensional and it is proven in~\cite{CH03,thesebenoit} that 
the Minkowskian norm 
of the linear form $H_{\Phi}$ does not depend on $\Phi$, again provided suitable conditions on the decay of $\Phi^* g - b$ 
at infinity hold. 
Namely, if one imposes a basis $\left(V_{(\mu)}\right)_{\mu=0 \cdots n}$ of $\mcK_0$ and defines $p_{(\mu)} = H_{\Phi}(V_{(\mu)})$,
then the quantity $\eta(p_{(\mu)},p_{(\mu)}) = p_{(0)}^2 - \sum_{i=1}^n p_{(i)}^2$ is independant of $\Phi$.
We can now define an important notion. We first observe that, if the hyperbolic metric is written as in~(\ref{b}), then the functions
$$V_{(0)} = \sqrt{1+r^2} \;,\; V_{(i)} = x^i$$
form a basis of $\mcK_0$, where the $x^i$ are the Cartesian coordinates of $\mathbb{R}^n$.
\begin{engdef}
Let $(M,g)$ an asymptotically hyperbolic manifold, with a spherical conformal boundary. Let $\Phi$ be a diffeomorphism 
$[R,\infty)\times N \longrightarrow M_{ext}$, with a radial coordinate $r$.  
The \emph{energy-momentum vector} of $g$ (or of $(M,g)$) is the vector $\mathbf{p_g}$ with components
$p_{(0)} = H_{\Phi}(\sqrt{1+r^2})$, $p_{(i)} = H_{\Phi}(x^i)$.
If $\mathbf{p_g}$ is timelike or null and future-directed, its \emph{mass} is the number $m_g$ such that 
$m_g = \sqrt{\eta(\mathbf{p_g},\mathbf{p_g})}$.
\end{engdef}
Using the Minkowskian form $\eta$ of signature $(+,-,\cdots,-)$, we can say that $\mathbf{p_g}$ is \emph{spacelike} 
(resp. \emph{null}, \emph{timelike}) if $\eta(\mathbf{p_g},\mathbf{p_g}) < 0$ (resp. $=0$, $>0$). 
Moreover, for a timelike energy-momentum vector, we say that $\mathbf{p_g}$
is \emph{future-directed} (resp. \emph{past-directed}) if the first coordinate $p_{(0)}$ is positive (resp. negative).
It is important to note that from the invariance property stated above, the $\eta$-norm of an energy-momentum $\mathbf{p_g}$ does not 
depend on the diffeomorphism $\Phi$, provided the suitable decay properties of $\Phi^* g - b$ stated in~\cite{CH03,thesebenoit}
are satisfied, although the vector $\mathbf{p_g}$ itself transforms under isometries of the hyperbolic space. 

Note that when one considers the initial data $(g,k)$ of an asymptotically anti-de Sitter spacetime, 
one can similarly define the corresponding energy-momentum vector, see e.g. the work of D. Maerten~\cite{thesedaniel, Mae06} 
for more details.

\medskip

A rather more straightforward way exists to define the notions of mass and of energy-momentum for a S.A.H. manifold $(M,g)$.

Under the definition (\ref{sah}) above, one introduces the notion of \emph{mass-aspect function} of $g$,
defined on the conformal infinity of $(M,g)$ (here $\mathbb S^{n-1}$) as the trace with respect to $h_0$ of the mass-aspect tensor $h$:
$$\mu_{h_0} := \mathrm{tr}_{h_0}(h) = h_0^{AB}h_{AB}\;.$$
One can now define the energy-momentum vector (see also \cite{Wan01}) $\mathbf{p_g}$ 
as the vector of $\mathbb R^{1,n}$:
$$\mathbf{p_g} := \left(\int_{\mathbb S^{n-1}}\mu_{h_0}(x)d\sigma_{h_0} , 
\int_{\mathbb S^{n-1}}\mu_{h_0}(x)\;x\; d\sigma_{h_0} \right)\;,$$
where $d\sigma_{h_0}$ is the volume associated to the round metric $h_0$ of $\mathbb S^{n-1}$. 
This definition coincides with the previous definition
of~\cite{CH03},  up to a constant positive factor, which does not affect the causal character of the vector
in the Minkowski spacetime.

This terminology of ``mass'' and ``energy-momentum'' is meaningful in mathematical general relativity (see~\cite{Hum10} for the mass)
and one of the main concerns of mathematical relativists is to prove that these quantities are in fact ``positive'' or 
``well-oriented'' under local geometric assumptions corresponding to the positivity of the density of energy~\cite{ADM61}.
In the present context of A.H. manifolds, the statement is: 
\begin{conj}
Let $(M,g)$ be a complete, $n$-dimensional Riemannian, asymptotically hyperbolic manifold whose conformal infinity is the
$(n-1)$-unit sphere. Assume that the scalar curvature $R_g$ of $g$ satisfies $R_g \geq -n(n-1)$.
Then the energy-momentum $\mathbf{p_g}$ is timelike future-directed, unless $(M,g)$ is isometric to the hyperbolic space.  
\end{conj}
This general result is yet an open question. However, the last decade has seen significant progress towards the proof of it.
Indeed, Wang (\cite{Wan01}) for S.A.H. manifolds, then Chru\'sciel and Herzlich (\cite{CH03}) for general asymptotically hyperbolic manifolds,
have independantly proved the following version of the ``positive energy-momentum theorem'':
\begin{Theorem}\label{WCH} Let $(M,g)$ be a complete Riemannian manifold, asymptotically hyperbolic and spin, with dimension $n \geq 3$, and whose
scalar curvature satisfies the inequality $R_g \geq -n(n-1)$. 
Then, the energy-momentum vector $\mathbf{p_g}$, if it exists and is non-zero,
is timelike future-directed. It is zero if and only if $(M,g)$ is isometric to the hyperbolic space $\mathbb H^n$.
\end{Theorem}

Without the spin assumption, only partial results are known; one of them follows from L. Andersson, M. Cai and G. Galloway,
as in~\cite{ACG08}:
\begin{Theorem}\label{ACG} Let $(M,g)$ be a complete Riemannian manifold, S.A.H., with dimension $n$ with $3 \leq n \leq 7$, 
whose scalar curvature
satisfies $R_g \geq -n(n-1)$. If moreover the mass-aspect function does not change of sign (strictly speaking) 
on the conformal boundary,
then this sign is positive, or zero if and only if the manifold $(M,g)$ is isometric to the hyperbolic space $\mathbb H^n$.
\end{Theorem}
This statement is far to be satisfying, especially when one compares the knowledge on the equivalent statement in the asymptotically
flat case. We can make however the following remark:
\begin{engrk} The assumption on the constant sign of the mass-aspect function in Theorem~\ref{ACG} implies that the energy-momentum vector is causal.

Indeed, for a mass-aspect function $u_0$ of (strict) constant sign, and for all $i$, one has:
$$ \left(\int x_i \frac{|u_0|d\sigma_{h_0}}{\int |u_0|d\sigma_{h_0}}\right)^2 \leq 
\int x_i^2 \frac{|u_0|d\sigma_{h_0}}{\int |u_0|d\sigma_{h_0}}\;,
$$
from the Jensen inequality applied to the measure $|u_0|d\sigma_{h_0}$ and to the convex function $X \mapsto X^2$.
The integrals are here computed on the conformal (spherical) infinity $\{x_1^2 + \cdots + x_n^2 = 1\}$, 
which yields
$$\sum_{i=1}^n \left(\int x_i u_0(x)d\sigma_{h_0}\right)^2 \leq \left(\int u_0(x) d\sigma_{h_0}\right)^2\;.$$
\end{engrk}

This result encourages us to seek functions $u_0$, non-radial, to play the role of mass-aspect functions of A.H. metrics with arbitrary
energy-momentum vectors. This is the object of the next section.

\section{Existence theorem}
\label{sectionthmexistence}

This section is devoted to establish the following result:

\begin{Theorem}
Let $\mathbf{p}$ be a vector of the Minkowski spacetime $\mathbb R^{1,n}$. Then there exists a compact set $K \subset \mathbb R^n$ 
and a Riemannian metric $g$, asymptotically hyperbolic, conformally flat defined on $\mathbb R^n\setminus K$, 
with constant scalar curvature $R_g = -n(n-1)$, such that the energy-momentum vector of $g$ exists, and coincides with $\mathbf{p}$.
\end{Theorem}

\proof We first notice that if $g$ is such a metric, then it takes the form $g=u^{\frac{4}{n-2}} b$ in a neighborhood of 
the boundary at infinity, where $b$ is the hyperbolic metric, hence $u$ solves the Yamabe equation
$$-\;4\;\frac{n-1}{n-2}\;\Delta_b u + R_b\; u = R_g\; u^{\frac{n+2}{n-2}}\;,
$$
where $\Delta_b$ is the Laplace-Beltrami operator of $b$.
Given the above conditions on the scalar curvature, $u$ has to solve the equation:
\begin{equation}\label{eqpsi2}\Delta_b u = \frac{n(n-2)}{4}\left(u^\frac{n+2}{n-2} -u\right)\;.\end{equation}
We introduce now coordinates $r$ (radial coordinate) and $\theta$ (latitude) such that the hyperbolic metric reads
$$b = \frac{dr^2}{1+r^2} + r^2 (d\theta^2 +\sin^2\theta \;\sigma_{n-2})\;,$$
where $\sigma_{n-2}$ is the standard metric of the $(n-2)$-dimensional unit sphere. 
The coordinate $r$ is related to the defining function $\rho$ by $r^{-1}= \sinh \rho$.
We seek $u$, solution of~(\ref{eqpsi2}), as a power series in $\frac1r$:
$$u = 1 + \sum_{k=0}^\infty \frac{u_k}{r^{n+k}}\;.$$
Our motivation to do so is that if the power series converges near the conformal boundary,
then the metric $g$ is S.A.H., with a mass-aspect function that coincides with $u_0$ up to a positive constant factor.

We start by looking for coefficients $u_k$ defined on the sphere $\Sphere^{n-1}$ that depend only on the latitude coordinate $\theta$.
In this case, the energy-momentum vector of $g$ takes the form:
\begin{equation}\label{EM}\mathbf{p_g} =
\lambda \left(\int_{\theta=0}^{\pi}u_0(\theta)\sin \theta d\theta\ , 0\ ,\cdots, 0\ , \int_{\theta=0}^{\pi}u_0(\theta)\sin \theta \cos \theta d\theta  \right)\;,
\end{equation}
where $\lambda$ is a positive constant.
From this, since $u$ depends only on the coordinates $r$ and $\theta$, one has:
$$\Delta_b u = (1 + r^2) \frac{\partial^2 u}{\partial r ^2} + \left(\frac{n-1 + n r^2}{r}\right)\frac{\partial u}{\partial r} +
 \frac{1}{r^2}\left(\frac{\partial^2 u}{\partial \theta ^2} + 
(n-2)\frac{\cos \theta}{\sin \theta}\frac{\partial u}{\partial \theta}\right)\;,$$
which yields, after a term-by-term differentiation in the series:
$$  \Delta_b u = \sum_{k = 0}^\infty w_k(\theta)r^{-(n+k)}\;,$$
where
\begin{equation}\label{wfu}w_k = (k+1)(k+n)u_k + k(k+n-2)u_{k-2} + u_{k-2}'' + (n-2)\frac{\cos \theta}{\sin \theta}u_{k-2}'\;,\end{equation}
for all $k \geq 0$, with the convention $u_{-1} = u_{-2} = 0$, and where the symbols $'$ and $''$ 
indicate first and second derivatives with respect to the variable $\theta$.

In order to handle the right-hand-side term of~(\ref{eqpsi2}), we write $u^\frac{n+2}{n-2}$ as
$$u^\frac{n+2}{n-2} = 1 + \sum_{p=1}^{+\infty}\left(\begin{array}{cc}\frac{n+2}{n-2}\\
 p
\end{array}\right) \left(\sum_{k=0}^\infty \frac{u_k}{r^{n+k}}\right)^p\;,
$$ 
where $\left(\begin{array}{cc}\frac{n+2}{n-2}\\ p \end{array}\right)$ is the combinatorial coefficient 
$$\left(\begin{array}{cc}\frac{n+2}{n-2}\\
 p
\end{array}\right) = \frac{\frac{n+2}{n-2}\cdots \left(\frac{n+2}{n-2} -1\right)\left(\frac{n+2}{n-2} -p +1\right)}{p!}\;.
$$
One can also write
$$u^\frac{n+2}{n-2} = 1 + \sum_{p=1}^{+\infty}\left(\begin{array}{cc}\frac{n+2}{n-2}\\
 p
\end{array}\right) \sum_{k=0}^{\infty}\left(\sum_{k_1+\cdots+k_p=k}u_{k_1}\cdots u_{k_p}\right) \frac{1}{r^{pn+k}}\;,
$$ 
We then identify the coefficients of the series $u^{\frac{n+2}{n-2}} = 1 + \sum_{l=0}^{\infty} \frac{v_l}{r^{n+l}}$.
If we denote $\mcE_{n,k}$ the set of elements $(p,l)$ such that $p \geq 2$, $l \geq 0$, $pn + l = n + k$, we have
$$v_k = \left(\frac{n+2}{n-2}\right)u_k + \sum_{(p,l) \in \mcE_{n,k}}\left(\begin{array}{cc}\frac{n+2}{n-2}\\ p \end{array}\right)
\sum_{l_1+\cdots+l_p = l}u_{l_1}\cdots u_{l_p}\;,
$$
and the equation~(\ref{eqpsi2}) together with the expression of $w_k$ yields, for all $k \geq 0$:
\begin{eqnarray}\label{wfu1}(k+1)(k+n) u_k + k(k+n-2) u_{k-2} + (n-2)\frac{\cos \theta}{\sin \theta}u_{k-2}' + u_{k-2}'' \hspace{2cm}\\ 
\hspace{3cm}\nonumber{= \frac{n(n-2)}{4}(v_k - u_k)\;.}
\end{eqnarray}
Note that when considering this equation for $k=0$, one gets that the coefficient $u_0(\theta)$ can be chosen freely, whereas for $k=1$,
the equation above forces $u_1(\theta)$ to be $0$.
Then, from the equation~(\ref{wfu1}) and the formula for $v_k$, all the coefficients $u_k$ are completely determined for $k \geq 2$
by the coefficients of lower rank.
More precisely, one can write
\begin{eqnarray}\label{wfu3}
\left((k+n+1)k + \frac{n(n-2)}{4}\right)u_k +  k(k+n-2) u_{k-2} + (n-2)\frac{\cos \theta}{\sin \theta}u_{k-2}' + u_{k-2}'' & \\
\nonumber{= \frac{n(n-2)}{4}P_k(u_0,\cdots,u_{k-1})\;,} &
\end{eqnarray}
where 
$$P_k = \sum_{p \geq 2\;,\; (p-1)n \leq k}\left(\begin{array}{cc}\frac{n+2}{n-2}\\ p \end{array}\right)
\sum_{l_1+\cdots+l_p = k - (p-1)n}u_{l_1}\cdots u_{l_p}\;. $$

In the sequel, we choose
$$u_0(\theta) = \beta + \cos \theta\;,$$
where $\beta \in \mathbb R$.
When doing this choice, we will obtain energy-momentum vectors given by the formula~(\ref{EM}) 
which are either timelike, or null, or even spacelike, depending on the values of $\beta$.
We use the subsequent fact, valid for all $\beta$:
\begin{lemm} With this choice of $u_0$, for all $k \geq 1$, $u_k(\theta)$ is a polynomial in the variable $\cos \theta$, 
with degree at most $k-1$.
\end{lemm}

\proof This property is trivial for $k=1$; let $k \geq 2$, such that $\deg(u_l)|_{\cos \theta}\leq l-1$, $\forall l \in \{1,\cdots,k-1\}$.
Then $\left(u_{k-2}'' + (n-2)\frac{\cos \theta}{\sin \theta}u_{k-2}' + k(k+n-2)u_{k-2}\right)$ is a polynomial in $\cos \theta$,
with degree $\leq k-1$ by assumption, and since $\deg(u_0)|_{\cos \theta}=1$.

Concerning the part $P_k = P_k(u_0,\cdots,u_{k-n})$, each term $\sum_{l_1 + \cdots l_p = l}u_{l_1} \cdots u_{l_p}$ (with $p \geq 2$) 
is made of products $u_{l_1} \cdots u_{l_p}$, which have, by assumption, a degree $\leq (l_1 +1) + \cdots + (l_p + 1) = l + p = 
k - (p-1)n + p = k - (n-1)(p-1) + 1$ which is less than $k-1$ as desired since $n \geq 3$ and $p \geq 2$. 
The formula $(\ref{wfu3})$ enables us to conclude the proof of the lemma.

\qed

\medskip

In the whole sequel, if $v = \sum_i v_i (\cos \theta)^i$ is a polynomial in $\cos \theta$, we note $|v|_1 = \sum_i |v_i|$.
\begin{lemm} There exists a number $\alpha > 0$ such that, for all positive integer $k$, one has
\begin{equation} |u_k|_1 \leq \frac{\alpha ^k}{(k+1)^2}\;.\end{equation}
\end{lemm}

\proof If $v = \sum_i v_i (\cos \theta)^i$ is a polynomial in $\cos \theta$, one has
$\frac{\cos \theta}{\sin \theta}v' = \sum_i -i v_i (\cos \theta)^i$
and $v'' = \sum_i v_i (-i (\cos \theta)^i + i(i-1) \sin^2 \theta (\cos \theta)^{i-2})$.

Thus one obtains the following bounds on the norms
$|.|_1$ of $\frac{\cos \theta}{\sin \theta} v'$ and of $v''$:
$$\left|\frac{\cos \theta}{\sin \theta}v'\right|_1 \leq d |v|_1 $$
and
$$|v''|_1 \leq d(d-1)|v|_1\;,$$
where $d$ is the degree of $v$ as a polynomial in $\cos \theta$.

Let us now assume that we have found a suitable $\alpha$, such that the statement of the lemma is satisfied at the order 
$l \in \{0,\cdots,k-1\}$
(it is obviously the case for $k=1$).
We infer from the result obtained above on $u_l$ that
\begin{eqnarray}\label{lemme2}\left|u_{k-2}'' +  (n-2)\frac{\cos \theta}{\sin \theta}u_{k-2}' + k(k+n-2)u_{k-2}\right|_1 \hspace{4cm} \\
\hspace{2cm}\nonumber{\leq \frac{\alpha^{k-2}}{(k-1)^2}((k-1)(k-2)+(n-2)(k-1)+k(k+n-2))\;.}
\end{eqnarray}
On the other hand, from our assumptions, one has
\begin{eqnarray}
|P_k|_1 & \leq & 
\sum_{(p,l) \in \mcE_{n,k}}\left|\left(\begin{array}{cc}\frac{n+2}{n-2}\\ p 
\end{array}\right)\right|
\sum_{l_1+\cdots+l_p = l} |u_{l_1}|_1\cdots |u_{l_p}|_1 \\
& \leq & 
\nonumber{\sum_{(p,l) \in \mcE_{n,k}}\left|\left(\begin{array}{cc}\frac{n+2}{n-2}\\ p 
\end{array}\right)\right|
\sum_{l_1+\cdots+l_p = l} \frac{\alpha^{l}}{(l_1+1)^2 \cdots (l_p + 1)^2}\;,}
\end{eqnarray}
since for all $(p,l) \in \mcE_{n,k}$, one has $l = k - (p-1)n \leq k-1$.
We now wish to evaluate the sums which appear above, for $(p,l) \in \mcE_{n,k}$, as
$$\mcS_p(l) := \sum_{l_1 + \cdots + l_p = l} \frac{1}{(l_1 + 1)^2} \cdots \frac{1}{(l_p + 1)^2}\;. 
$$
We first need the following general result, simply obtained by a decomposition into simple elements of a rationnal fraction (see in Appendix):
\begin{equation}\label{s2}
 \sum_{r=0}^q \frac{1}{(r+1)^2 (q-r+1)^2} \leq \frac{\pi^2}{(q+2)^2}\;.
\end{equation}
From this, there is no difficulty to find an upper bound for $\mcS_p(l)$ as:
$$\mcS_p(l) \leq \frac{\pi^{2(p-1)}}{(l+p)^2}\;.
$$
Hence, one can write
$$|P_k|_1 \leq \sum_{p\geq 2\;,\; (p-1)n \leq k} \left|\left(\begin{array}{cc}\frac{n+2}{n-2}\\ p 
\end{array}\right)\right| \frac{\pi^{2(p-1)} \alpha^{k-(p-1)n}}{(k-(p-1)n + p)^2}\;,
$$
so we need to estimate the combinatorial term $ \left|\left(\begin{array}{cc}\frac{n+2}{n-2}\\ p 
\end{array}\right)\right| $. We can in fact show (see Appendix) that it bounded by $C e^p$, where $C$ is a positive constant.

On the other hand, we wish to find an upper bound for the ratio $\left(\frac{k+1}{k-(p-1)n+ p}\right)^2$. We find (see Appendix) that
$n^2$ is a valid upper bound of this ratio for all $p \geq 2$ such that $(p-1)n \leq k$.

Combining these results, we can write
$$|P_k|_1 \leq
\frac{\alpha ^k}{(k+1)^2}\sum_{p \geq 2\;,\; (p-1)n \leq k} C_n \left(\frac{e \pi^2}{\alpha^n}\right)^{p-1}\;,
$$ 
where $C_n$ is a number that depends only on $n$. Thus, for all $\alpha$ such that $e \pi^2 < \alpha^n$,
one concludes that
\begin{equation}\label{pk}
 |P_k|_1 \leq \left(\frac{C_n}{1-\frac{e \pi^2}{\alpha^n}} \frac{e \pi^2}{\alpha^n}\right) \frac{\alpha^k}{(k+1)^2}\;,
\end{equation}
and the term inside the bracket can be made as small as desired for $\alpha$ large enough, for any given $n$.

Putting together~(\ref{wfu3}),~(\ref{lemme2}) and~(\ref{pk}), we obtain
\begin{eqnarray}\label{uk}
|u_k|_1 \leq 
\left[\frac{(k-1)(k+n-4) + k(k+n-2)}{k^2 + (n+1)k + \frac{n(n-2)}{4}}\times \frac{1}{\alpha^2} \times \left(\frac{k+1}{k-1}\right)^2 
\right. \hspace{2cm}
\\
\hspace{4cm}\nonumber{\left. + \frac{\frac{n(n-2)}{4}}{k^2 + (n+1)k + \frac{n(n-2)}{4}} \times \frac{C_n}{1-\frac{e \pi^2}{\alpha^n}} \frac{e \pi^2}{\alpha^n}\right]
\frac{\alpha^k}{(k+1)^2} \;.}
\end{eqnarray}
In particular, the term inside the bracket admits an upper bound independant of $k$ and can be made as small as desired when choosing 
a sufficiently large value of $\alpha$.  
Therefore, if we choose $\alpha$ large enough (independantly of $k$) we have
\begin{equation} 
|u_k|_1 \leq \frac{\alpha^k}{(k+1)^2}\;,
\end{equation}
so that we can conclude by induction on $k$ the proof of the lemma.

\qed

\medskip

We infer from this result that $|u_k|_{\infty} :=\sup_{\theta} |u_k(\theta)| \leq \frac{\alpha^k}{(k+1)^2}$ for all $k$, 
since the $u_k$ are polynomials in $\cos \theta$.
We have therefore obtained that $u (r,\theta)$, solution of $(\ref{eqpsi2})$, 
as a power series in $1/r$, with a convergence radius less than $\alpha^{-1}$; 
in other words, we have constructed a metric $g=u^{\frac{4}{n-2}} b$, defined on a neighborhood of the conformal (spherical) boundary,
of the form
$$M_{ext}= (\alpha,+\infty)\times \mathbb{S}^{n-1}\;,$$
with the desired properties, and non-radial.

For $\beta = 0$, from the formula (\ref{EM}), the energy-momentum vector $\mathbf{p_g}$ is spacelike, whereas it is timelike
for large enough values of $|\beta|$.
Moreover, these computations can be conducted also for functions $u_0$ of the form $u_0(\theta) = \beta + \gamma \cos \theta$, 
with $\gamma \in \mathbb R$, leading to a similar conclusion.
Hence, when $\beta$ and $\gamma$ vary in $\R$, the energy-momentum $\mathbf{p_g}$ of the resulting metrics describe the two dimensional subspace
$\R \times \{0\}^{n-1} \times \R$ of the Minkowski space $\R^{1,n}$.
One can now reach any vector $\mathbf p$ of $\mathbb R^{1,n}$ when using furthermore the action of isometries ``at infinity'' (of the hyperbolic
space), namely using the group $O(1,n)$.

\qed

\section{Applications}\label{appli}
\subsection{Optimality of a positive energy-momentum theorem with boundary}

In this section, we explain why the result obtained in the previous section illustrates that the Theorem 4.7 in~\cite{CH03} 
is somehow optimal.

Let us first recall the statement of this theorem:
\begin{Theorem}[Chru\'sciel-Herzlich 03]
Let $(M,g)$ be a $n$-dimensional Riemannian manifold, complete, spin, where $g$ is a $\mathcal C^2$, and $M$ has
a (inner) compact, non-empty boundary of mean curvature
$$\Theta \leq n-1\;,$$
and where the scalar curvature of metric $g$ satisfies
$$R_g \geq -n(n-1)\;.$$
Under suitable assumptions on the asymptotic behaviour and if the conformal infinity is the $(n-1)$-unit sphere, then the 
energy-momentum vector $\mathbf{p_g}$ is timelike future-directed.
\end{Theorem}
As we would see, this result is no longer valid if the assumption on the mean curvature $\Theta$ of the inner boundary is removed.

The first obvious counter-example is given by the hyperbolic space $\mathbb H ^n$ itself, as the energy-momentum vector 
is trivial in this case, while the scalar curvature is $R_b = -n(n-1)$. Indeed for the hyperbolic metric, 
one can check that there is no compact hypersurface in $\mathbb H ^n$ with mean curvature $\Theta \leq n-1$. 
Indeed, let $\Sigma$ be a closed hypersurface in $\mathbb H ^n$.
Consider any $(n-1)$-sphere that contains $\Sigma$ in its interior and tangent to $\Sigma$ at at least one point $p$.
Then, the value of the mean curvature of $\Sigma$ at $p$ is at least the value of the mean curvature of this $(n-1)$-sphere, 
which is precisely $(n-1) \coth r$ for a sphere of radius $r$, hence stricly larger than $n-1$.

More generally speaking, let us consider a S.A.H. $n$-dimensional manifold $(M,g)$, endowed with a defining function $\rho$.
Outside some compact of $M$, the metric takes the expression, in the limit $\rho \rightarrow 0$:
$$g = \sinh ^{-2} \rho \left(d \rho^2 + h_0 + \frac{\rho^n}{n} h + o (\rho^n)\right)\;,$$
where $h_0$ is the standard metric on the $(n-1)$-unit sphere.
Let us now compute the mean curvature of ``large spheres'' in $M$, whose center lies at the
origin of the coordinate system of the chart at infinity considered here.

Defining the coordinate $s$ by $\sinh s = \sinh ^{-1} \rho$, the S.A.H. metric takes the following form as $s \rightarrow \infty$:
$$g = ds^2 + \sinh^2 s \ h_s\;,$$
where $h_s = h_0 + \frac{2^n}{n}  e^{-n s} h + o (e^{-n s})$.
The second fundamental form of the sphere $S_s$ of euclidean radius $s$ centered at the origine  of this coordinate system reads:
$$I\!I_{S_s} = \frac12 \frac{\partial}{\partial s}\left(\sinh^2 s\ h_s\right)
= 2 \sinh s \cosh s\ h_s + \sinh^2 s \frac{\partial h_s}{\partial s}\;,$$
and the mean curvature is therefore
$$\Theta _{S_s} = \left(\sinh^2 s \ h_s\right)^{-1} I\!I_{S_s} = (n-1)\frac{\cosh s}{\sinh s}
+ \frac12 h^{-1}_s \frac{\partial h_s}{\partial s}\;.
$$
Then, one has $\partial_s h_s = - 2^n e^{-ns} h + o(e^{-ns})$ as $s \rightarrow +\infty$. Thus, one has
$$\Theta _{S_s} = (n-1)\frac{\cosh s}{\sinh s} - 2^{n-1} e^{-ns} \mu_{h_0} + o(e^{-ns})\;,$$
and we recover the mass-aspect function $\mu_{h_0}$ which does not depend on $s$. 
Hence, for $n \geq 3$, one obtains the asymptotic behaviour
of the mean curvature of large $(n-1)$-spheres:
$$
\Theta _{S_s} = (n-1)\left(1 + 2 e^{-2s}\right) + O(e^{-3s})\;.
$$
This last quantity is strictly larger than $n-1$ for all $s$ large enough. 
One notices that the above expansion is valid for all S.A.H. metrics,
no matter of the assumptions on the scalar curvature or on the mass-aspect function $\mu_{h_0}$.

Hence, if one replaces the condition of the mean curvature of the inner boundary in Chru\'sciel-Herzlich's Theorem
by a condition
$$\Theta \leq \alpha (n-1)$$
for $\alpha  > 1$, then the condition is completely general: given such an number $\alpha$, 
every S.A.H. metric has an inner compact, boundaryless hypersurface whose mean curvature is less than $\alpha (n-1)$,
in particular this is the case for metrics constructed in the previous section, thus with arbitrary energy-momentum vector.

\subsection{A gluing result for asymptotically hyperbolic metrics}
We show here that the gluing results of P.T. Chru\'sciel and E. Delay in~\cite{CD09} apply when one, instead of using the family of \emph{boosted}
Kottler (Schwarzschild-anti de Sitter) metrics, uses the family constructed in section~\ref{sectionthmexistence} as models 
for the asymptotic region. 

\begin{figure}
 \begin{center}
\hspace{-2cm}{
 \psfrag{M}{\Huge  $\!\!M$}
 \psfrag{Mun}{\Huge$ \!\!M_1$}
 \psfrag{Mdeux}{\Huge$ M_2$}
 \psfrag{gtilde}{\Huge$ \tilde g$}
 \psfrag{g}{\Huge $\!g$}
 \psfrag{gp}{\Huge $\!\!\!\!\!\mathring g_Q$}
\hspace{2.5cm}\resizebox{3in}{!}{\includegraphics{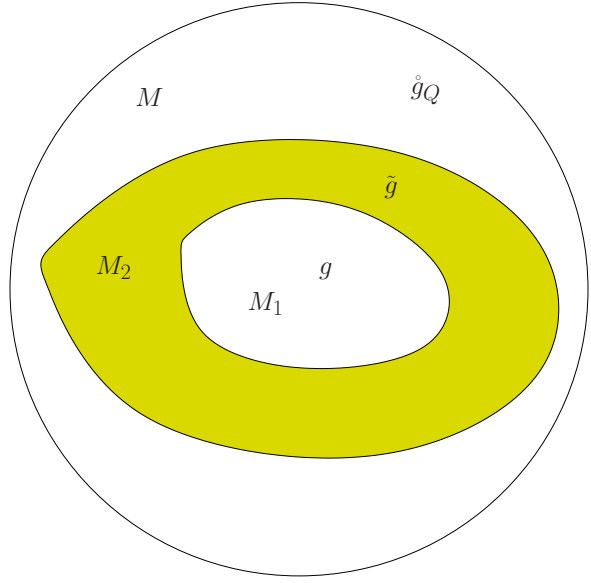}}
}
\caption{The set-up for the gluing construction.
\label{massAHgluing}}
\end{center}
\end{figure}

Let us give a quick overview of the Corvino-Schoen's gluing principle, at least for \emph{static} initial data.
We start by considering a non-compact Riemannian manifold $(M,g)$ which asymptotes to a reference metric $b$ in the asymptotic region.
(We will assume for simplicity that there is only one asymptotic end). The metric $g$ is furthermore assumed to have constant scalar curvature
$R_g = R_b$. The gluing result consists in finding compact regions $M_1 \subset \subset M_2 \subset M$, 
a family of ``model'' metrics $(\mathring g_Q)_{Q \in \mcF}$ in the asymptotic region of $M$ and a new metric $\tilde g$ on $M$
satisfying the following properties, as illustrated in figure~\ref{massAHgluing}:
\begin{itemize}
 \item $\tilde g$ coincides with $g$ on $M_1$,
 \item $\tilde g$ coincides with $\mathring g_Q$ on $M\setminus M_2$ for some parameter $Q$ in $\mcF$,
 \item the scalar curvature $R_{\tilde g}$ of $\tilde g$ is constant on $M$ and equal to $R_b$.
\end{itemize}

The proof of such results consists in interpolating $g$ and $\mathring g _Q$ by defining $g_Q = (1-\chi)g + \chi \mathring g _Q$,
with a troncature function $0 \leq \chi \leq 1$ which is equal to $0$ in $M_1$ and to $1$ in $M\setminus M_2$, and then
to find a compactly supported (in $M_2\setminus M_1$) perturbation $\delta g$ such that we recover $R(g_Q + \delta g) = R_b$.
The kernel $\mcK_0$ of the adjoint of the linearized scalar curvature operator plays here an essential obstruction role, and only
a big enough family $\mcF$ of parameters for an suitable family $(\mathring g _Q)$ of models at infinity ensures the existence of a
suitable $Q \in \mcF$ which solves the perturbation problem. In fact (see~\cite{CS06,CD03}), the dimension of $\mcF$ has to be at least equal
to the dimension of $\mcK_0$.    

Going back to the situation of this paper, the metrics exhibited in section~\ref{sectionthmexistence} form
an ``admissible'' family (in the sense of\cite{CS06}) of models for the infinity, in the sense that they satisfy the asymptotic decay requirements imposed in~\cite{CD09}, 
and their global charges (in fact their energy-momentum vectors) describe a non-trivial $(n+1)$-dimensional open set in the $(n+1)$-Minkowski 
space of all the possible global charges. In fact, from what we saw in the previous section, the energy-momentum vectors of 
our family reach the whole Minkowski space, which means that the space of parameters $\mcF$ coincides with $\mathbb{R}^{1,n}$.

\medskip

We start by introducing some notations and definitions that appear in~\cite{CD09}. If $\rho$ is a defining function
of the (spherical) conformal infinity, we define 
$$M_{\varepsilon} := \{x \in M\;,\; \rho(x) < \varepsilon\}$$
and the annulus
$$A_{\delta,\varepsilon} := M_{\varepsilon}\setminus M_{\delta}$$
for $0 < \delta < \varepsilon$. We also denote by $\nabla$ the Levi-Civit\`a connection of the (hyperbolic) metric $b$.

For every vector $p_{(\mu)}$ in the Minkowski space $\mathbb{R}^{1,n}$, we define $\mathring g _{p_{(\mu)}}$ to be the metric that we have 
constructed in section~\ref{sectionthmexistence} with an energy-momentum vector $p_{(\mu)}$. 
For all $n \geq 5$, we define 
$$\alpha_n := \max\left(8, \frac{8+n}{2}\right)\;.$$
We can now state the result, directly adapted from Theorem 1.2 of~\cite{CD09}:
\begin{Theorem}
 Let $n \geq 5$, $l \in \mathbb N$ with $l > \lfloor \frac{n}{2} \rfloor + 4$, $\lambda \in  (0,1)$ and let $\alpha > \alpha_n$.
Let $g$ be a $C^{l,\lambda}$-asymptotically hyperbolic metric on $M$ with constant scalar curvature $-n(n-1)$ and an energy-momentum vector 
$p^0_{(\mu)}$ such that the conditions hold:
$$|g - \mathring g_{p^0_{(\mu)}}|_b + \cdots + |\nabla^{l-2}(g - \mathring g_{p^0_{(\mu)}})|_b = O(\rho^{\alpha})\;,
$$
and
$$|g - \mathring g_{p^0_{(\mu)}}|_b + \cdots + |\nabla^{l}(g - \mathring g_{p^0_{(\mu)}})|_b = O(1)\;.
$$
Then, there exists $\delta_0 > 0$ such that for all $\delta \in (0,\delta_0]$, there exists a metric $\tilde g$ on $M$ satisfying the following
conditions:
\begin{itemize}
 \item $\tilde g$ has constant scalar curvature $-n(n-1)$,
 \item $\tilde g$ coincides with $g$ on $M\setminus M_{4\delta}$,
 \item $\tilde g$ coincides with the metric $\mathring g_{p_{(\mu)}}$ on $M_{\delta}$.
\end{itemize}
Moreover, $\tilde g$ is smooth if $g$ is.

\end{Theorem}

\proof 
Most of the details of the proof can be found in~\cite{CD09}. In particular, the results of the section 3 in that paper hold.
For section 4, we just replace the family of boosted Kottler metrics by the family of metrics constructed here in the previous section, and 
$g$ is here allowed to have any energy-momentum vector (not necessarily timelike). 
\qed

\begin{engrk} A few comments can be done: \end{engrk}
 \begin{itemize}
  \item For $n \geq 3$, one may also find a adapted version of the Theorem 1.1 of~\cite{CD09}, replacing Kottler metrics by metrics 
$\mathring g_{p_{(\mu)}}$ constructed above.
  \item  As in~\cite{CD09} with the boosted Kottler metrics, one can obtain a metric $\mathring g_{p_{(\mu)}}$ of energy-momentum 
$p_{(\mu)}$ as close as desired from the energy-momentum $p^0 _{(\mu)}$ of $g$, provided the number $\delta_0$ is small enough.
This means that one can deform any asymptotically hyperbolic metric $g$ with $R_g = -n(n-1)$ into a metric $\tilde g$, which is also
asymptotically hyperbolic (and even S.A.H.), with $R_{\tilde g} = -n(n-1)$, which is conformally flat in a neighborhood of 
the conformal infinity, and whose energy-momentum vector is as close as desired from the energy-momentum vector of $g$.
  \item As pointed out in~\cite{CD09}, the asymptotic decay rate of $g$ towards $\mathring g_{p^0_{(\mu)}}$ required 
for the theorem to hold ($\alpha > \alpha_n$) is undesirably high, especially in low dimensions, compared to the decay rate needed to define the 
energy-momentum vector (where $\alpha > n/2$ suffices).
 \end{itemize}
%
%
Note that the construction performed in this work may well be extended to mass-aspect functions which take the form of polynomials 
in $\cos \theta$
and $\sin \theta$ or even more generally to analytic functions with conditions on the coefficient growth, although the simpler construction
presently considered already suffices to achieve our goals concerning the energy-momentum vector. It would however be of 
interest to study this further and find whether there are conditions for a function to be realized as a mass-aspect function of some 
S.A.H. metric, satisfying or not the scalar curvature requirement $R_g = -n(n-1)$.

\appendix
\section{Appendix}
\begin{Lemma}
For all $k \geq 0$, $n \geq 3$ and $p \geq 2$ such that $(p-1)n \leq k$, the following inequality holds:
$$\left(\frac{k+1}{k-(p-1)n + p}\right)^2 \leq n^2\;.
$$
\end{Lemma}

\proof The quantity on the left-hand side is positive and also takes the form $\left(1 - \frac{(p-1)(n-1)}{k+1}\right)^{-2}$, 
and by assumption, $(p-1)(n-1) \leq k\frac{n-1}{n}$, hence we have 
$$\left(\frac{k+1}{k-(p-1)n + p}\right)^2 \leq \left(1 - \frac{k}{k+1}\frac{n-1}{n}\right)^{-2} = 
\left(\frac{n+k}{(k+1)n}\right)^{-2}\;,
$$
and the result follows.
\qed

\begin{Lemma}
There exists a positive constant $c$ such that, for all positive integers $p$ and $n$, one has
$$ \left|\left(\begin{array}{cc}\frac{n+2}{n-2}\\ p 
\end{array}\right)\right| \leq c e^p\;.$$
\end{Lemma}

\proof
We first note that $\left|\left(\begin{array}{cc}\frac{n+2}{n-2}\\ p 
\end{array}\right)\right|$ is less than the quantity $\frac{1}{p!}
\max \left\{\left(\frac{n+2}{n-2}\right)^p\;,\; p^p \right\}$.
The result then comes from the Stirling formula $p ! \sim p^p e^{-p} \sqrt{2 \pi p}$.

\qed

\begin{Lemma}
For all $q \geq 0$,
$$\sum_{r=0}^q \frac{1}{(r+1)^2 (q-r+1)^2}  \leq \frac{\pi^2}{(q+2)^2}\;.
$$
\end{Lemma}

\proof
Let us denote by $S_2(q)$ the left-hand side term.
The result follows from the decomposition:
\begin{eqnarray}\frac{1}{(r+1)^2 (q-r+1)^2} = \hspace{9.5cm}\\
\hspace{1cm}\nonumber{ \frac{1}{(q+2)^2}\left(\frac{1}{(r+1)^2} + \frac{1}{(q-r+1)^2}\right)  + 
\frac{2}{(q+2)^3}\left(\frac{1}{r+1} + \frac{1}{q-r+1}\right)\;. }
\end{eqnarray}
Indeed, one can write 
$$S_2(q) = \frac{2}{(q+2)^2}\sum_{r=0}^q \left(\frac{1}{(r+1)^2} + \frac{2}{(q+2)(r+1)}\right)
\leq \frac{6}{(q+2)^2}\sum_{r=0}^q \frac{1}{(r+1)^2}\;,
$$
and the last term above is bounded by $\frac{6 \zeta (2)}{(q+2)^2} = \frac{\pi^2}{(q+2)^2}$.
\qed

\bibliographystyle{amsplain}

\begin{thebibliography}{11}

\bibitem{ACG08}
L.~Andersson, M.~Cai, and G.J. Galloway, \emph{Rigidity and positivity of mass
  for asymptotically hyperbolic manifolds}, Ann. Henri Poincar\'e \textbf{9}
  (2008), 1--33.

\bibitem{ADM61}
R.~Arnowitt, S.~Deser, and C.~Misner, \emph{Coordinate invariance and energy
  expressions in general relativity}, Phys. Rev. \textbf{122} (1961),
  997--1006.

\bibitem{Bar86}
R.~Bartnik, \emph{The mass of an asymptotically flat manifold}, Commun. Pure
  Appl. Math. \textbf{39} (1986), 661--693.

\bibitem{BI04}
R.~Bartnik and J.~Isenberg, \emph{The {C}onstraint {E}quations},  (2004),
  arXiv: gr-qc/0405092v1.

\bibitem{BC97}
R.~Beig and P.T. Chru\'sciel, \emph{{Killing Initial Data}}, Class. Quantum
  Grav. \textbf{14} (1997), A83--A92.

\bibitem{CD03}
P.T. Chru\'sciel and E.~Delay, \emph{On mapping properties of the general
  relativistic constraints operator in weighted function spaces, with
  applications}, M\'em. Soc. Math. Fr. (2003), no.~94.

\bibitem{CD09}
\bysame, \emph{Gluing constructions for asymptotically hyperbolic manifolds
  with constant scalar curvature}, Comm. Anal. Geom. \textbf{17} (2009), no.~2,
  343--381.

\bibitem{CH03}
P.T. Chru\'sciel and M.~Herzlich, \emph{The mass of asymptotically hyperbolic
  {R}iemannian manifolds}, Pacific J. Math \textbf{212} (2003), no.~2,
  231--264.

\bibitem{Cor00}
J.~Corvino, \emph{Scalar curvature deformation and a gluing construction for
  the {E}instein constraint equations}, Comm. Math. Phys. \textbf{214} (2000),
  no.~1, 137--189.

\bibitem{CS06}
J.~Corvino and R.M. Schoen, \emph{On the asymptotics for the vacuum {E}instein
  constraint equations}, J. Diff. Geom. \textbf{73} (2006), no.~2, 185--217.

\bibitem{theseromain}
R.~Gicquaud, \emph{\'{E}tude de quelques probl\`emes d'analyse et de
  g\'eom\'etrie sur les vari\'et\'es asymptotiquement hyperboliques}, Ph.D.
  thesis, Universit\'e Montpellier 2, 2009.

\bibitem{Hum10}
E.~Humbert, \emph{Relativit\'e g\'en\'erale (d'apr\`es {M}. {V}augon) et
  quelques probl\`emes math\'ematiques qui en sont issus},  (2010), arXiv:
  1004.2402v4 [math.DG].

\bibitem{thesedaniel}
D.~Maerten, \emph{Aspects math\'ematiques du moment-\'energie en relativit\'e
  g\'en\'erale}, Ph.D. thesis, Universit\'e Montpellier 2, 2006.

\bibitem{Mae06}
\bysame, \emph{Positive energy-momentum theorem in asymptotically anti-de
  {S}itter space-times}, Ann. Henri Poincar\'e \textbf{7} (2006), 975--1011.

\bibitem{thesebenoit}
B.~Michel, \emph{Invariants asymptotiques en g\'eom\'etrie conforme et
  g\'eom\'etrie {CR}}, Ph.D. thesis, Universit\'e Montpellier 2, 2010.

\bibitem{Mon76}
V.~Moncrief, \emph{Space-time symmetries and linearization stability of the
  {E}instein equations. {II}}, J. Math. Phys. \textbf{17} (1976), no.~10,
  1893--1902.

\bibitem{ST07}
Y.~Shi and L-F.~Tam, \emph{Asymptotically hyperbolic metrics on the unit ball 
with horizons}, Manuscripta Math. \textbf{122} (2007), 97--117.

\bibitem{Smontecatini}
R.M. Schoen, \emph{Variational theory for the total scalar curvature functional
  for {R}iemannian metrics and related topics}, Topics in Calculus of
  Variations (Montecatini Terme, 1987), Lecture Notes in Mathematics, vol.
  1365, Springer Berlin / Heidelberg, 1989, pp.~120--154.

\bibitem{Wan01}
X.~Wang, \emph{The mass of asymptotically hyperbolic manifolds}, J. Diff. Geom
  \textbf{57} (2001), no.~2, 273--279.

\end{thebibliography}

\end{document}